\documentclass[12pt]{article}

\usepackage{pstricks,pst-node,graphicx}
\usepackage{amssymb,amsmath}

\newcommand{\be}{\begin{equation}}
\newcommand{\ee}{\end{equation}}
\newcommand{\eea}{\end{eqnarray}}
\newcommand{\bea}{\begin{eqnarray}}
\newcommand{\xin}{\xi_N}
\newcommand{\szn}{Z_{\xin}}
\newcommand{\zn}{Z_{\xin}(u,v,w)}
\newcommand{\zns}{ \sum\limits_{\xin} \zn}

\newcommand{\snzn}{ {1 \over {2^N}}\sum\limits_{\xin} \zn}

\newcommand  \nt {{\noindent}}
\newcounter{secnum}[section]
\begin{document}


\title{{{\bf Coloured Hard-Dimers}}
\footnotetext{\newline{\bf Key words and phrases}: Coloured hard-dimers, combinatorics, generating function
\newline
{\it Mathematics Subject Classification}: 05A15, 60C05}}

\author{{\bf Maria Simonetta Bernabei} \\ and \\
{\bf Horst Thaler} \\[1ex] Department of Mathematics and Informatics,  \\
University of Camerino, \\
Via Madonna delle Carceri 9,
I--62032, Camerino (MC), Italy;\\
{\small simona.bernabei@unicam.it, horst.thaler@unicam.it}}

\date{}

\maketitle

\abstract{{An averaged generating function for coloured hard-dimers is being investigated by proving estimates for the latter. Furthermore, two different enumerating problems and their distributions are studied numerically.

}}

\section{Introduction}
 \setcounter{secnum}{\value{section}
 \setcounter{equation}{0}
 \renewcommand{\theequation}{\mbox{\arabic{secnum}.\arabic{equation}}}}
\newtheorem{remark}{Remark}[section]
\newtheorem{lemma}{Lemma}[section]
\newtheorem{proposition}{Proposition}[section]
\newtheorem{theorem}{Theorem}[section]

Coloured hard-dimers have recently appeared in the context of causally triangulated (2+1)-dimensional gravity. Using special triangulations of spacetime, it was shown in \cite{BeLoZa} that the generating function of the one step propagator can as well be expressed through the generating functions of coloured hard-dimer configurations.

In this article we go in a slightly different direction and study the behaviour of an averaged generating function for coloured hard-dimers. In this manner we touch some combinatorial problems, which so far have not been tackled.

The paper is  organized as follows. In section 2, we define the notion of coloured hard-dimers, give the statement of the problem and prove estimates for the averaged generating function. In section 3 additional material for numerics around coloured hard-dimers is presented.

\section{ Description of the model and results}
 \setcounter{secnum}{\value{section}
 \setcounter{equation}{0}
 \renewcommand{\theequation}{\mbox{\arabic{secnum}.\arabic{equation}}}}

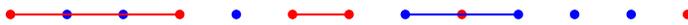
\begin{figure}[h]
\begin{center}
\psset{xunit=1.5cm,yunit=0.75cm}
\begin{pspicture}(0,0)(7,2)
\psdots[dotstyle=*,linecolor=red](0.5,0.5)(2.0,0.5)(3,0.5)(3.5,0.5)(4.5,0.5)(6.5,0.5)
\psdots[dotstyle=*,linecolor=blue](1,0.5)(1.5,0.5)(2.5,0.5)(4,0.5)(5,0.5)(5.5,0.5)(6,0.5)
\psline[linecolor=red](0.5,0.5)(2,0.5)
\psline[linecolor=red](3,0.5)(3.5,0.5)
\psline[linecolor=blue](4,0.5)(5,0.5)
\end{pspicture}
\caption{\label{hdimer} A hard-dimer, $N=13$}
\end{center}
\end{figure}

Let $\xi_N$ be a sequence of length $N$ of blue and red sites on the one-dimensional lattice
$\mathbb{Z}$. Let us define the notion of a coloured hard-dimer on $\xi_N$: a coloured hard-dimer is a sequence of red or blue dimers, which cannot intersect each other (``hardness" property). A dimer in turn is an edge linking two nearest sites of the same colour. In Fig.\ref{hdimer} an example of a coloured hard-dimer is given.

For any configuration $\xin$, the generating function associated with it is of the form
\be
\label{nr1}
Z_{\xin}(u,v,w)=\sum_{D} u^{n_b(D)} v^{n_r(D)} w^{n_{br}(D)},\quad u,v,w\in (0,\infty),
\ee
where $D$ is an hard-dimer on $\xin$, $n_b(D)$ and $n_r(D)$ denote the number
of blue and red dimers respectively on  $D$, and $n_{br}(D)$ the number of crossings between dimers and sites of different colour. Moreover, let us denote by $\gamma_b(D)$ and $\gamma_r(D)$ the number of blue and red sites of $D$ respectively, which do not belong to any dimer. We call such kind of sites ``single points". Naturally the following constraint
\be
\label{nr2}
2 n_b(D) + 2n_r(D) + n_{br}(D)+  \gamma_b(D) +\gamma_r(D)=N
\ee
holds.
It is useful to define further variables $t$ and $s$, where $t$ corresponds to the number of sites occupied by dimers, and $s$ to the number of dimers:
\be
\label{nr3}
t= N-  \gamma_b(D) - \gamma_r(D),
\qquad
s= n_b(D) + n_r(D)
\ee
In the above example $n_b(D)=1,\,  n_r(D)=2,\,  n_{br}(D)=3,\,  \gamma_b(D)=3$ and $\gamma_r(D)=1$.

For any $N$, we want to prove estimates from above and from below, for the mean of the generating functions $Z_{\xin}$ defined in (\ref{nr1}). In achieving this we calculate explicitly, for a given $N$, the sum of all the generating functions $\szn$ over $\xin$, rescaled by the factor $2^N$, i.e. the number of all the sequences $\xin$. We get so an ``averaged" generating function
$${1 \over {2^N}} \sum_{\xin} \zn.$$
The main result of this section is the following:

\begin{theorem}

For any $N$, the following estimates

\be
\label{nr5}
1 - {{C_1} \over {2^N \sqrt{N}}}+{C \over {\sqrt{N}}} \left ( {{A+1}\over 2}\right )^N  \leq
\snzn \leq
1 - {{1} \over { \sqrt{N}}}+ {1 \over {\sqrt{N}}} \left ( {A+1} \right )^N
\ee
hold, where $A \equiv \sqrt{u+v} +{w \over 2}$ and $C_1, C$ are constants not depending on $N$, given
as
$$C_1=e^{-\frac{1}{6}}\sqrt{\frac{2}{\pi}},\; {C}=\frac{C_1}{2} \left ( {1} \wedge {{ w+2} \over { \sqrt{u+v} +{\frac w2} +1}} \right )$$

\end{theorem}

\begin{remark}

Note that $A >0$, so that the upper bound in ($\ref{nr5}$) is always exponentially divergent, as $N \to \infty$, whereas the lower bound increases exponentially  for $A>1$ and large $N$, i.e. for the class of parameters $u,v,w$, such that $\sqrt{u+v}+{\frac w2}>1$, otherwise it tends to $1$.

\end{remark}

In order to prove the main result, we first need a Lemma, where we get an explicit formula for the average ${1 \over{2^N}} \sum_{\xin} Z_{\xin}$, for any $N$. We are able to get an exact formula
by summing over the variables $t$ and $s$ defined above in ($\ref{nr3}$).

\begin{lemma}

The average ${1 \over{2^N}} \sum_{\xin} Z_{\xin}$ has the following
explicit expression, for any $N$:

$$
\snzn =$$
\be
\label{nr6}
1+  \sum\limits_{t=1}^N \sum\limits_{s=1}^{[{\frac t2}]} {N-t+s \choose s} {t-s-1 \choose s-1} \left ({{u+v}\over 4}\right )^s \left( {w \over 2} \right )^{t-2s}
\ee
where $[\cdot]$ denotes the integer part.
\end{lemma}

\nt
{\bf Proof of Lemma 2.1:} Consider a  hard-dimer $D$ of length $N$, and set $n_b(D)=k$, $n_r(D)=h$, $n_{br}(D)=m$, $\gamma_b(D)=\gamma_b$ and $\gamma_r(D)=\gamma_r$.
By  ($\ref{nr2}$) and  ($\ref{nr3}$) the following equalities
\be
\label{nr7}
2k+2h+m+\gamma_b+\gamma_r=N, \qquad  t=2k+2h+m
\ee
hold.

First we fix the number of blue and red dimers and that one of single points $\gamma_b$ and $\gamma_r$, without assigning any length   to dimers, even though, by ($\ref{nr3}$), the variable $t$ on $D$ is done. In this manner the dimers of the same colour are considered indistinguishable. Then we calculate all possible permutations of $k$ blue dimers, $h$ red dimers, $\gamma_b$ blue and $\gamma_r$ red single points, i.e.

$$
{{(k+h+\gamma_b+\gamma_r)!}\over{k! \;  h! \; \gamma_b! \gamma_r!}}
$$

Now, one has to distinguish two cases: $k=h=0$ and $k+h \neq 0$. If $k=h=0$, then $m=0$, because this means that there are no dimers. Therefore $\gamma_b+ \gamma_r =N$ and the combinatorics is
$$
 \sum\limits_{\gamma_b=0}^{N} {N \choose{ \gamma_b}} = 2^N
 $$

In the other case ($k+h \neq 0$) we assign, for any given permutation of coloured dimers and single points, all the admissible lengths to each dimer, taking into account the fact that $t$ is given. The number of all these combinations is then
$$
{m+k+h-1 \choose k+h-1}
$$
Therefore we have
$$
\zns =$$
\be
\label{nr8}
\sum_{k,h,m: \atop 2k+2h+m=1}^N  \sum_{ \gamma_b + \gamma_r=N-2k-2h-m} {{(k+h+\gamma_b+\gamma_r)!}\over{k! \;  h! \; \gamma_b! \gamma_r!}} {m+k+h-1 \choose k+h-1} u^k v^h w^m
\ee
From ($\ref{nr7}$) we have that $k+h+\gamma_b+\gamma_r=N-k-h-m$ and $\gamma_r=N-2k-2h-m-\gamma_b$ and hence
$$
\zns =   $$
$$ \sum\limits_{k,h,m: \atop 2k+2h+m=1}^N  \sum\limits_{\gamma_b=0}^{N-2k-2h-m} {{(N-k-h-m)!}\over{k! \; h! \; \gamma_b! (N-2k-2h-m-\gamma_b)!}} {m+k+h-1 \choose k+h-1} u^k v^h w^m = $$
$$
\sum\limits_{k,h,m: \atop 2k+2h+m=1}^N  {{(N-k-h-m)!}\over{k! \; h! (N-2k-2h-m)!}} {m+k+h-1 \choose k+h-1} u^k v^h w^m \times $$
\be
\label{nr9}
 \sum\limits_{\gamma_b=0}^{N-2k-2h-m}  {N-2k-2h-m \choose \gamma_b}
\ee
In ($\ref{nr9}$) we have multiplied and divided the generic term of the sum by $(N-2k-2h-m)!$, so that the binomial formula  appears $$\sum_{\gamma_b=0}^{N-2k-2h-m} {N-2k-2h-m \choose \gamma_b}= 2^{N-2k-2h-m}$$
 and $\sum_{\xin} {\szn}$ becomes

\be
\label{nr10}
2^N \left [  1+    \sum\limits_{k,h,m: \atop 2k+2h+m=1}^N   {{(N-k-h-m)!}\over{k! h! (N-2k-2h-m)!}} {m+k+h-1 \choose k+h-1} \left ({u\over 4}\right )^k \left ({v\over 4}\right )^h \left ({w\over 2}\right )^m \right ]
\ee
Note that the scaling factor $2^N$ naturally appears in ($\ref{nr10}$). Now, by multiplying and dividing the generic term of the sum in ($\ref{nr10}$) by $(k+h)!$ we get
$$
 \snzn = $$
$$  1 + \sum\limits_{k,h,m: \atop 2k+2h+m=1}^N   {N-k-h-m \choose  k+h} {m+k+h-1 \choose k+h-1} {k+h \choose k}  \left ({u\over 4}\right )^k \left ({v\over 4}\right )^h \left ({w\over 2}\right )^m  $$
 We perform the variable changements $s=k+h$ and $t=2k+2h+m$, where $s$ and $t$ were already defined. We get thus
$$
 \snzn = $$
 $$ 1 +   \sum\limits_{t=1}^N   \sum\limits_{s=1}^{[{\frac t2}]}  {N-t+s \choose  s} {t-s-1 \choose s-1}
\left ({w\over 2}\right )^{t-2s}  \sum\limits_{k=0}^s {s \choose k} \left ({u \over 4}\right )^k \left ({v\over 4}\right )^{s-k}   =
$$
\be
\label{nr11} 1 + \sum\limits_{t=1}^N   \sum\limits_{s=1}^{[{\frac t2}]}  {N-t+s \choose  s} {t-s-1 \choose s-1}   \left ({u+v\over 4}\right )^s  \left ({w\over 2}\right )^{t-2s}
\ee
In the last sum in ($\ref{nr11}$) we used again the binomial formula
\be \sum_{k=0}^{s} {s \choose k}  \left({u\over 4}\right)^k \left({v\over 4}\right )^{s-k} =
\left ({{u+v}\over 4}\right )^s \nonumber 
\ee
 Therefore we get formula ($\ref{nr6}$) where only the indices $s$ and $t$ appear. The Lemma is so proved.
\hfill{$\Box$}
\medskip

\noindent
{\bf Proof of Theorem 2.1:}  First we write the sum on the right hand side of formula ($\ref{nr11}$) in a more suitable form:
\be
\label{nr12}
{N-t+s \choose s} {t-s-1 \choose s-1}={{s {2s \choose s}}\over{N{N-1 \choose t-s-1}}}
{N \choose 2s} {N-2s \choose t-2s}, \quad s\geq 1,\, t \geq 1
\ee

In fact, it holds
$$ {N-t+s \choose s} {t-s-1 \choose s-1} = $$
$$ {{(N-t+s)!} \over  {s! (N-t)!}} \cdot
 {{(t-s-1)!} \over  {(s-1)! (t-2s)!}} \cdot
  {{N!} \over  {N!}} \cdot  {{(N-2s)!} \over  {(N-2s)!}} \cdot
  {{(2s)!} \over  {(2s)!}} =
$$
$$ {{
   {{(2s)! } \over  {s! (s-1)!}}
   }\over
   { {{N! } \over  {(N-t+s)! (t-s-1)!}}
   }}
{N \choose 2s} {N-2s \choose t-2s} =
$$
$$ {{
s {2s \choose s}}
\over { N {N-1 \choose t-s-1}}}
{N \choose 2s} {N-2s \choose t-2s}
$$
Now we estimate the binomial coefficients ${2s \choose s}$ and ${N-1 \choose t-s-1}$ from above and from below, more precisely, first we prove that
\be
\label{nr13}
{{C_1 2^{2s}} \over {2^N \sqrt{N}}} \leq
 {{
s {2s \choose s}}
\over { N {N-1 \choose t-s-1}}}
\leq {{2^{2s}} \over {\sqrt{N}}}
\ee
where $C_1 = e^{-{\frac 16}}\sqrt{{2 \over \pi}}$.
In order to estimate the binomial coefficient ${2s \choose s}$ we use Stirling's formula
$$
s!=  s^{s} e^{-s} \sqrt{2 \pi s}   \; e^{\lambda_s},
$$
where
$$
{1 \over {12s+1}} \leq \lambda_s \leq  {1 \over {12s}}
$$
We may write
\be
\label{nr14}
{2s \choose s}=  {{2^{2s}} \over {\sqrt{\pi s}}} \; e^{\lambda_{2s} - 2 \lambda_s}
\ee
Taking into account the following inequalities for the errors
$$- {1 \over {6s}} \leq  \lambda_{2s} - 2 \lambda_s \leq 0$$
we get estimates from above and from below for $s {2s \choose s}$
\be
\label{nr15a}
{{e^{-{\frac 16}}} \over {\sqrt{\pi}}} \; 2^{2s} \leq
s {2s \choose s} \leq 2^{2s} \sqrt{s}
\ee

For the binomial coefficient ${N-1 \choose t-s-1}$ we use the fact that its maximum is assumed at $t-s-1 = \left [{{N-2} \over 2} \right ]$. Hence
$$
1 \leq {N-1 \choose t-s-1} \leq  {N-1 \choose \left [{{N-2} \over 2} \right ] }
$$
We can assume for simplicity  that $N$ is even. The case of $N$ odd can be treated with little more effort.
Then, as before, we use Stirling's formula
$$
 N {N-1 \choose {\frac N2}-1} = {{N!}\over{\left (  {\frac N2}-1\right )! \left (  {\frac N2}\right )! }} = $$
$$  {\frac N2}  {{N!}\over{ \left (  {\frac N2}\right )!^2 }} \leq  {\frac N2} {{2^N}\over{\sqrt{\frac N2}}} =
2^N \sqrt{{\frac N2} }
$$
Thus we obtain
$$N \leq N {N-1 \choose t-s-1} \leq  N {N-1 \choose \left [{{N-2} \over 2} \right ] } \leq 2^N \sqrt{\frac{N}{2}},
$$
that yields with  ($\ref{nr15a}$) the formula ($\ref{nr13}$).
Therefore we are able to estimate the mean ${1 \over {2^N}} \sum_{\xin} Z_{\xin}$ from below and above using the same series
$$
1 +  {{C_1} \over {2^N \sqrt{N}}}\sum\limits_{t=1}^N \sum\limits_{s=1}^{[{\frac t2}]} {N \choose 2s} {N-2s \choose t-2s} \left (u+v \right )^s \left( {w \over 2} \right )^{t-2s} \leq  $$
 $$ \snzn \leq $$
\be
\label{nr15} 1 + {1 \over {\sqrt{N}}}\sum\limits_{t=1}^N \sum\limits_{s=1}^{[{\frac t2}]} {N \choose 2s} {N-2s \choose t-2s} \left (u+v \right )^s \left( {w \over 2} \right )^{t-2s}
\ee
It remains to evaluate the series
\be
\label{nr16}
\sum\limits_{t=1}^N \sum\limits_{s=1}^{[{\frac t2}]} {N \choose 2s} {N-2s \choose t-2s} \left (u+v \right )^s \left( {w \over 2} \right )^{t-2s}
\ee
We then state the next Lemma

\begin{lemma}

The following equality

$$ \sum\limits_{t=1}^N \sum\limits_{s=1}^{[{\frac t2}]} {N \choose 2s} {N-2s \choose t-2s} \left (u+v \right )^s \left( {w \over 2} \right )^{t-2s}= $$
\be
\label{nr17}
 {\frac 12} \left [ \left (  \sqrt{u+v} + {w \over 2} + 1 \right )^N + \left ( - \sqrt{u+v} + {w \over 2} + 1 \right )^N  \right ] -1
\ee
holds, for any $N$.
\end{lemma}

\nt
{\bf Proof of Lemma 2.2:}  Consider the sum ($\ref{nr16}$) for $s$ and $t$ starting from $0$ and denote it by $\ast$. Note that in $\ast$ there is an extra term at $t=s=0$ which equals 1. In this case we can apply the following combinatorial equality (see \cite{Fe})
\be
\label{nr18}
\sum\limits_\nu {n \choose \nu} {n-\nu \choose k-\nu} t^\nu = {n \choose k} (t+1)^k,
\ee
where the sum is over any $\nu$ such that $0 \leq \nu \leq k \leq n$.
The summation with respect to $s$ in $\ast$ coincides
with the sum ($\ref{nr18}$) restricted to the even indices $\nu$. Owing to this fact we rewrite the sum over $s$ in a more suitable form:

$$ \sum\limits_{s=0}^{[{\frac t2}]} {N \choose 2s} {N-2s \choose t-2s} \left (\sqrt{u+v} \right )^{2s} \left( {w \over 2} \right )^{t-2s}= $$
$$   \sum\limits_{s=0}^t {{1+(-1)^s}\over 2}  {N \choose s} {N-s \choose t-s} \left (\sqrt{u+v} \right )^s \left( {w \over 2} \right )^{t-s} = $$
$$ {\frac 12}  \sum\limits_{s=0}^t {N \choose s} {N-s \choose t-s} \left ( {{2 \sqrt{u+v}} \over w}  \right )^s \left( {w \over 2} \right )^{t} + $$
$$ {\frac 12} \sum\limits_{s=0}^t  {N \choose s} {N-s \choose t-s} \left (-{{2 \sqrt{u+v}} \over w} \right )^s \left( {w \over 2} \right )^{t}
$$
We apply now formula ($\ref{nr18}$)
$$ \sum\limits_{s=0}^{[{\frac t2}]} {N \choose 2s} {N-2s \choose t-2s} \left (\sqrt{u+v} \right )^{2s} \left( {w \over 2} \right )^{t-2s}= $$
$$
{\frac 12}{N \choose t} \left [ \left (  {{2\sqrt{u+v}}\over w} + 1 \right )^t  + \left ( - {{2\sqrt{u+v}}\over w} + 1 \right )^t   \right ]
 \left( {w \over 2} \right )^{t} =$$
 \be
 \label{nr19}
   {\frac 12}{N \choose t}  \left [
 \left (  {\sqrt{u+v}} +  {w \over 2}  \right )^t
  + \left (  - {\sqrt{u+v}} +  {w \over 2} \right )^t
  \right ]
\ee
Summing the last term in ($\ref{nr19}$) over $t$ we get again a binomial formula
$$ {\frac 12} \sum\limits_{t=0}^N {N \choose t}  \left [
 \left (   {\sqrt{u+v}} + {w \over 2}  \right )^t
  + \left ( - {\sqrt{u+v}} + {w \over 2}  \right )^t
  \right ] = $$
  $$  {\frac 12}    \left [
 \left (  {\sqrt{u+v}}  +{w \over 2} + 1\right )^N
  + \left (  - {\sqrt{u+v}} +{w \over 2} + 1\right )^N
  \right ]
$$
The Lemma 2.2. is so proved.
\hfill{$\Box$}

\medskip
Now we are able to finish the proof of  the Theorem
\medskip

\nt
{\bf Upper bound:} Taking into account ($\ref{nr15}$) and ($\ref{nr17}$) we obtain an upper bound for ${1 \over{2^N}} \sum_{\xin} Z_{\xin}$
$$
\snzn \leq$$
$$ 1 - {1 \over {\sqrt{N}}}+ {1 \over {2 \sqrt{N}}}  \left (   {\sqrt{u+v}} +{w \over 2} + 1\right )^N
\left [  1 +  \left ( {{   - {\sqrt{u+v}} +{w \over 2} +1 } \over {  {\sqrt{u+v}}+{w \over 2} +1}}\right )^N \right ]
$$
Since
$$ {{  - \sqrt{u+v} +{w \over 2}+1 } \over {  \sqrt{u+v} +{w \over 2}+1}}\leq 1,
$$
we have proved the upper estimate in ($\ref{nr5}$).

\nt
{\bf Lower bound:} Analogously, we obtain as lower bound for ${1 \over{2^N}} \sum_{\xin} Z_{\xin}$
$$  1 - {{C_1} \over {2^N \sqrt{N}}}+  {{C_1} \over {2^{N+1}\sqrt{N}}}
\left [ \left (    {\sqrt{u+v}} +{w \over 2} + 1\right )^N +  \left (  - {\sqrt{u+v}} +{w \over 2} + 1\right )^N \right ]= $$
$$ 1 - {{C_1} \over {2^N \sqrt{N}}}+ {{C_1} \over {2 \sqrt{N}}}  \left [ \left ( {{\sqrt{u+v}} \over 2}    +{w \over 4} + {\frac 12} \right )^N +
 \left (  - {{\sqrt{u+v}} \over 2}  +{w \over 4} + {\frac 12} \right )^N \right ] =$$
\be
\label{nr20}
1 - {{C_1} \over {2^N \sqrt{N}}}+ {{C_1} \over {2 \sqrt{N}}}  \left ({{\sqrt{u+v}} \over 2}  +{w \over 4} +  {\frac 12}  \right )^N
\left [ 1+ \left ( {{ - {{\sqrt{u+v}}}  +{w \over 2} +1 } \over {
  {{\sqrt{u+v}} +{w \over 2} +1} }}  \right )^N  \right ]\ee
 Therefore we get the lower bound ($\ref{nr5}$), with the constant
$$ {C}= {{C_1} \over 2} \quad \textnormal{for $- \sqrt{u+v} + {\frac w2} +1 \geq 0$} $$
$${C}={{C_1} \over 2} \left ( 1+ {{- {{\sqrt{u+v}} +{w \over 2} +1}} \over {
  {{\sqrt{u+v}}}   +{w \over 2} +1 }} \right )  =  \frac{C_1}{2} \left ( {{w+2} \over {
  {{\sqrt{u+v}}  +{w \over 2} + 1} }} \right ),
$$
for $- \sqrt{u+v}  + {\frac w2} + 1 < 0$. The Theorem is so proved.
\hfill{$\Box$}

\section{Numerics for coloured hard-dimers}
 \setcounter{secnum}{\value{section}
 \setcounter{equation}{0}
 \renewcommand{\theequation}{\mbox{\arabic{secnum}.\arabic{equation}}}}

In this section we provide some supplementary material regarding hard-dimer combinatorics. We consider two enumeration problems which are different from that in section 2. Numerical techniques are being employed and the results may indicate interesting directions for further study. In the first part we count the number of hard-dimers for each configuration $\xi_N$, $N$ fixed. The so obtained statistics is then shown in histograms.

In order to explain the algorithm we consider a particular configuration as shown in Fig.\ref{edimers}. \\[1ex]
\begin{figure}[h]
\begin{center}
\psset{xunit=1.5cm,yunit=0.6cm}
\begin{pspicture}(0,0)(5,4)
\psdots[dotstyle=*,linecolor=red](0,0.5)(1.5,0.5)(2.5,0.5)(3,0.5)(3.5,0.5)(4.5,0.5)
\psdots[dotstyle=*,linecolor=blue](0.5,0.5)(1,0.5)(2,0.5)(4,0.5)
\psline[linecolor=red]{|-|}(0,1)(1.5,1)
\psline[linecolor=blue]{|-|}(0.5,1.5)(1,1.5)
\psline[linecolor=blue]{|-|}(1,2)(2,2)
\psline[linecolor=red]{|-|}(1.5,2.5)(2.5,2.5)
\psline[linecolor=blue]{|-|}(2,3)(4,3)
\psline[linecolor=red]{|-|}(2.5,3.5)(3,3.5)
\psline[linecolor=red]{|-|}(3,4)(3.5,4)
\psline[linecolor=red]{|-|}(3.5,4.5)(4.5,4.5)
\end{pspicture}\hspace{2ex}
\parbox[b]{3.5cm}{\small Dimers for a particular configuration of sites are indicated as closed intervals.\\[1ex]}
\end{center}

\caption{\label{edimers} Dimers}
\end{figure}

Recall that a hard-dimer on a sequence of red and blue sites is a sequence of dimers which are non-overlapping. In order to count the number of hard-dimers on the fixed sequence $\xi_N$ of red and blue sites, one can proceed as follows:\\[1ex]
{\it Step 1.} Generate a list of all dimers, i.e. codify their colour as well as their starting and end points. E.g. for the dimers in Fig.\ref{edimers} this is achieved by defining the set
$$ED:=$$
$$\{\{1,4,r\},\{2,3,b\},\{3,5,b\},\{4,6,r\},\{5,9,b\},\{6,7,r\},\{7,8,r\},\{8,10,r\}\}$$
{\it Step 2.} For each $l\in \{1,\ldots,|ED|\}$ let $K(l)$ be the set of subsets of $ED$ containing $l$ elements. Do the following for $l=1$ to $l=|ED|$. Scan the elements of $K(l)$. If for an element in $K(l)$ the dimers are non-overlapping, increase the number of hard-dimers by one, otherwise keep it unchanged. At the end one has the number of hard-dimers for a particular configuration $\xi_N$. \\[0.5ex]

For fixed $N$, let $S_N$ be the set of all $\xi_N$ and let $\xi_N(D)$ denote the number of hard-dimers for $\xi_N$. This number is obtained following the steps 1-2 above. Then let $X=(\xi_N(D))_{\xi_N \in S_N}$ be the sample obtained from this numeric counting procedure. The histograms for the standardized sample $Y=(X-\rm{mean}(X))/{\rm st. deviation}(X)$ are depicted in Fig.\ref{hist2a} and Fig.\ref{hist2b} below.

\begin{figure}[ht]
\begin{center}
\includegraphics[scale = 0.75]{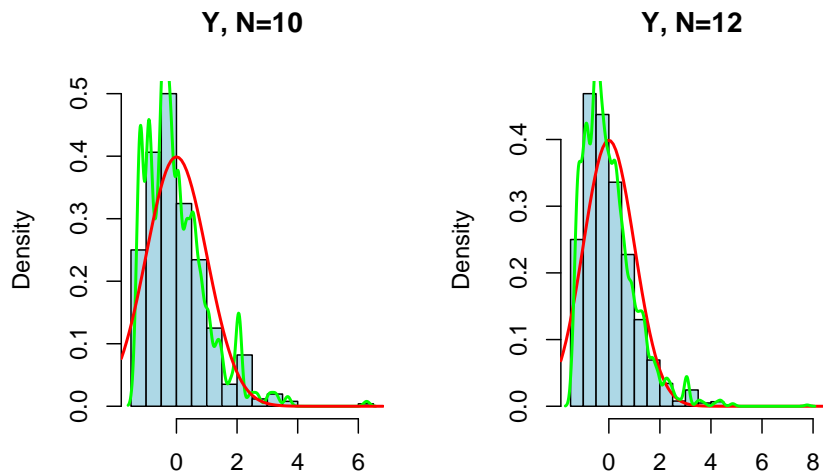}
\caption{\label{hist2a} Histograms of $Y$}
\end{center}
\end{figure}

\begin{figure}[ht]
\begin{center}
\includegraphics[scale = 0.7]{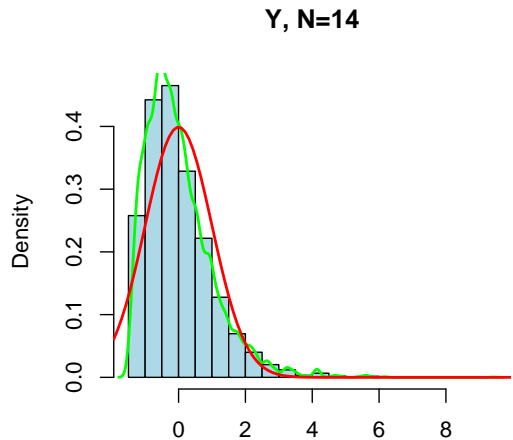}
\caption{\label{hist2b} Histogram of $Y$}
\end{center}
\end{figure}

\begin{remark}{
The red curves in the histograms of $Y$ above represent the standard normal distribution. On the other hand, the green curves are obtained by smoothing of $Y$ with a Gaussian kernel having standard deviation 0.1. With increasing $N$ the green curves become less oscillating and a tendency towards a definite distribution becomes visible. Unfortunately, the algorithm becomes very time consuming for $N$ larger 14 and it is not clear whether the distributions tend to a standard normal distribution. }
\end{remark}

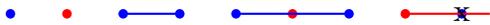
\begin{figure}[ht]
\psset{xunit=1.5cm,yunit=0.75cm}
\begin{pspicture}(0,0)(5.5,1)
\psdots[dotstyle=*,linecolor=red](1.0,0.5)(3,0.5)(4,0.5)
\psdots[dotstyle=*,linecolor=blue](0.5,0.5)(1.5,0.5)(2,0.5)(2.5,0.5)(3.5,0.5)(4.5,0.5)
\psline[linecolor=blue](2.5,0.5)(3.5,0.5)
\psline[linecolor=blue](1.5,0.5)(2,0.5)
\psline[linecolor=red](4,0.5)(4.75,0.5)\rput(4.5,0.5){x}
\end{pspicture}
\caption{\label{harddimero} A hard-dimer being omitted}
\end{figure}

In this second part we generate hard-dimers randomly. Roughly speaking, this is achieved by first choosing $N$ red or blue sites randomly, each with probability one half and then for each such configuration dimers are generated. Here again randomness enters in that one dimer and the next one are separated by some random choice. The algorithm can be implemented in the following steps:\\[1ex]
{\it Step 1.} Generate a sequence of red and blue sites. The colours appear with probability one half. \\[1ex]
{\it Step 2.} For the first site, depending on the colour, we choose randomly the base elements with attributes equal or right, see Fig.\ref{rdimer}. For example, if the colour is blue, br or be are picked at random with probability one half. 
\newpage
\begin{figure}[ht]
\begin{center}
\psset{xunit=1.5,yunit=0.8}
\begin{pspicture}(0,0)(5,5)
\psframe(0,0)(4.5,4.5)
\psdots[dotstyle=*,linecolor=red](1,0.75)\rput(1.75,0.75){re}
\psdots[dotstyle=*,linecolor=red](1,1.75)
\psline[linecolor=red](1,1.75)(1.25,1.75)\rput(1.75,1.75){rr}
\psline[linecolor=red](1,2.75)(0.75,2.75)
\psdots[dotstyle=*,linecolor=red](1,2.75)\rput(1.75,2.75){rl}
\psline[linecolor=blue](0.75,3.75)(1.25,3.75)\rput(1.75,3.75){rm}
\psdots[dotstyle=*,linecolor=red](1,3.75)
\psdots[dotstyle=*,linecolor=blue](3,0.75)\rput(3.75,0.75){be}
\psdots[dotstyle=*,linecolor=blue](3,1.75)
\psline[linecolor=blue](3,1.75)(3.25,1.75)\rput(3.75,1.75){br}
\psline[linecolor=blue](3,2.75)(2.75,2.75)
\psdots[dotstyle=*,linecolor=blue](3,2.75)\rput(3.75,2.75){bl}
\psline[linecolor=red](2.75,3.75)(3.25,3.75)\rput(3.75,3.75){bm}
\psdots[dotstyle=*,linecolor=blue](3,3.75)
\end{pspicture} \hspace{1ex}
\parbox[b]{4.5cm}{\small These are the building elements from which the random hard-dimers are constructed. \\ The acronyms mean: re $\doteq$ red-equal, rr $\doteq$ red-right, rl $\doteq$ red-left, rm $\doteq$ red-mixed, and similarly for blue.}
\caption{\label{rdimer} Base elements}
\end{center}
\end{figure}
{\it Step 3.} Suppose we have generated the configuration up to the $k^{\rm th}$ site and want to configure the $(k+1)^{\rm th}$ site. Here we have to distinguish the four possible combinations of the two colours. Suppose the colours are both blue. Then, given the element at the $k^{\rm th}$ site we have to choose the base element of the $(k+1)^{\rm th}$ site according to the following table: \\[2ex]
\begin{tabular}{|c|c|} \hline
$k^{\rm th}$ site is blue & $(k+1)^{\rm th}$ site is blue \\ \hline \hline
br & bl \\ \hline
bm & bm \\ \hline
be & random choice $\in \{\mbox{be,br}\}$  \\ \hline
bl & random choice $\in \{\mbox{be,br}\}$  \\ \hline
\end{tabular} \\[2ex]
In the case of a blue and a red site we take: \\[2ex]
\begin{tabular}{|c|c|} \hline
$k^{\rm th}$ site is blue & $(k+1)^{\rm th}$ site is red \\ \hline \hline
br & rm \\ \hline
bm & rd \\ \hline
be & random choice $\in \{\mbox{re,rr}\}$  \\ \hline
bl & random choice $\in \{\mbox{re,rr}\}$  \\ \hline
\end{tabular} \\[2ex]
The random choice happens with probability one half. Similarly for the other two cases, taking into account the properties of dimers. \\[1ex]
{\it Step 4.} For the last site the elements are chosen from the base elements with attributes left and equal only. If this is not possible, then a $\times$ is being inserted. In this way we indicate that the dimer remains unfinished and will therefore not be counted in the last step 5, see Fig.\ref{harddimero}. Concretely, we follow the rules: \\[2ex]
\begin{tabular}{|c|c|} \hline
penultimate site is blue & last site is blue \\ \hline \hline
br & bl \\ \hline
bm & x \\ \hline
be & be  \\ \hline
bl & be  \\ \hline
\end{tabular} \\[2ex]
and \\[2ex]
\begin{tabular}{|c|c|} \hline
penultimate site is blue & last site is red \\ \hline \hline
br & x \\ \hline
bm & rl \\ \hline
be & re  \\ \hline
bl & re  \\ \hline
\end{tabular} \\[2ex]
Similarly for the other two cases.\\[1ex]
{\it Step 5.} Count the number of dimers.

\begin{figure}[ht]
\begin{center}
\includegraphics[scale = 0.7]{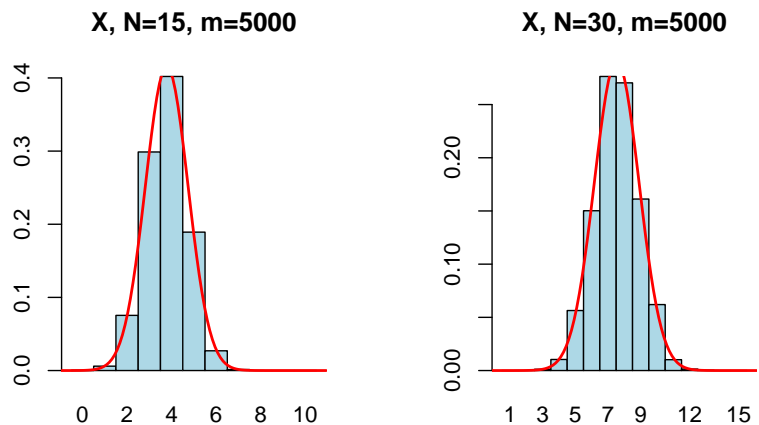}
\caption{\label{hist3} Histograms of $X$}
\end{center}
\end{figure}

Now, let $S^r_N$ denote a random sample of a certain size $|S^r_N|=m$ which is obtained by repeating the above procedure, steps 1-4, $m$ times. Again, $N$ is given by the number of sites. For $\eta\in S^r_N$ let $\eta(d)$ denote the number of dimers counted in $\eta$. Then in Fig.\ref{hist3} the histograms are shown for two samples $X=(\eta(d))_{\eta\in S^r_N}$.

\newpage

\begin{remark}
In the histograms for $X$ above, the red curves have been chosen to have the same mean and standard deviation as the corresponding histograms. The very good approximation by the histograms makes a central limit property visible.
\end{remark}


\begin{thebibliography} {9}



\bibitem{BeLoZa} Benedetti, D., Loll, R., Zamponi, F.: $(2+1)$-dimensional quantum gravity as the continuum limit of causal dynamical triangulations.  Phys. Rev. D  76  (2007),  no. 10, 104022, 26 pp.

\bibitem{Fe} Feller, W.: An introduction to probability theory and its applications. Vol. I. Third edition John Wiley $\&$ Sons, Inc., New York-London-Sydney 1968

\end{thebibliography}
\end{document}